\documentclass{article}%
\usepackage{amsmath}
\usepackage{amsfonts}
\usepackage{amssymb}
\usepackage{graphicx}%
\providecommand{\U}[1]{\protect\rule{.1in}{.1in}}
\newenvironment{proof}[1][Proof]{\noindent\textbf{#1.} }{\ \rule{0.5em}{0.5em}}
\begin{document}

\title{Completeness of eigenfunctions of Sturm-Liouville problems with
discontinuities at three points }
\author{{\Large Erdo\u{g}an \c{S}en}}
\date{}
\maketitle

{\scriptsize {Department of Mathematics, Faculty of Science and Letters, Nam\i
k Kemal University, 59030 Tekirda\u{g}, Turkey.}}

{\scriptsize e-mail: esen@nku.edu.tr}



ABSTRACT.

In this work, we study discontinuous Sturm-Liouville type problems with
eigenparameter dependent boundary condition and transmission conditions at
three interior points. A self-adjoint linear operator $A$ is defined in a
suitable Hilbert space $H$ such that the eigenvalues of such a problem
coincide with those of $A$. We show that the eigenvalues of the problem are
analytically simple, and the eigenfunctions of $A$ are complete in $H$.

\section{Introduction}

It is well-known that many topics in mathematical physics require the
investigation of eigenvalues and eigenfunctions of Sturm-Liouville type
boundary value problems. In recent years, more and more researches are
interested in the discontinuous Sturm-Liouville problems with
eigenparameter-dependent boundary conditions (see $[1-4])$. The literature on
this subject is voluminous and we refer to $[5-10]$. Various physics
applications of this kind problem are found in many literatures, including
some boundary value problem with transmission conditions that arise in the
theory of heat and mass transfer (see $[2,10,11]$). The study of the structure
of the solution in the matching region of the layer with the basis solution in
the plate leads to consideration of an eigenvalue problem for a second order
differential operator with piecewise continuous coefficients and transmission
conditions $\left[  12\right]  $.

Sturm-Liouville problems with transmission conditions have been studied by
many authors (see $[3,4,8,13]$). Adjoint and self-adjoint boundary value
problems with interface conditions have been stutied in $\left[  14,15\right]
$. Such problems with point interactions are also studied in $\left[
16\right]  $.

In this study, we also deal with the class of problems $\left(  1\right)
-\left(  9\right)  $, by means of a combination of the methods $\left[
5\right]  $,$\left[  13\right]  $,$\left[  17\right]  $ and$\left[  18\right]
$. In Section 1, a self-adjoint linear operator $A$ is defined in a suitable
Hilbert space $H$ such that the eigenvalues of the problem $\left(  1\right)
-\left(  9\right)  $ coincide with those of $A$. In Section 2, we prove that
the eigenvalues of the problem $\left(  1\right)  -\left(  9\right)  $ are
analytically simple. In Section 3, we prove that the eigenfunctions of $A$ are
complete in $H$. Note that each eigenfunction of the original problem and a
real number.

In this study, we consider a discontinuous eigenvalue problem which consists
of Sturm-Liouville equation
\begin{equation}
\tau u:=\left(  -p\left(  x\right)  u^{^{\prime}}\left(  x\right)  \right)
^{^{\prime}}+q(x)u\left(  x\right)  =\lambda u\left(  x\right)  \tag*{(1)}%
\end{equation}
on $I=\left[  -1,h_{1}\right)  \cup\left(  h_{1},h_{2}\right)  \cup\left(
h_{2},h_{3}\right)  \cup\left(  h_{3},1\right]  ,$ where $p\left(  x\right)
=p_{1}^{2}$ for $x\in\left[  -1,h_{1}\right)  $; $p\left(  x\right)
=p_{2}^{2},$ for $x\in\left(  h_{1},h_{2}\right)  ,$ $p\left(  x\right)
=p_{3}^{2},$ for $x\in\left(  h_{2},h_{3}\right)  $ and $p\left(  x\right)
=p_{4}^{2},$ for $x\in\left(  h_{3},1\right]  ;$ $p_{1},p_{2},p_{3},p_{4}$ are
nonzero real constants, $q\left(  x\right)  \in L^{1}\left(  I,%
%TCIMACRO{\U{211d} }%
%BeginExpansion
\mathbb{R}
%EndExpansion
\right)  $ and $\lambda\in%
%TCIMACRO{\U{2102} }%
%BeginExpansion
\mathbb{C}
%EndExpansion
$ is the eigenparameter; with the boundary condition%
\begin{equation}
\tau_{1}u:=\alpha_{1}u\left(  -1\right)  +\alpha_{2}u^{^{\prime}}(-1)=0,
\tag*{(2)}%
\end{equation}
the eigenparameter-dependent boundary condition%
\begin{equation}
\tau_{2}u:=\lambda\left[  \beta_{1}^{^{\prime}}u\left(  1\right)  -\beta
_{2}^{^{\prime}}u^{^{\prime}}(1)\right]  +\left[  \beta_{1}u\left(  1\right)
-\beta_{2}u^{^{\prime}}(1)\right]  =0, \tag{3}%
\end{equation}
and the transmission conditions%

\begin{align}
\tau_{3}u  &  :=u\left(  h_{1}+0\right)  -\alpha_{3}u\left(  h_{1}-0\right)
-\beta_{3}u^{^{\prime}}\left(  h_{1}-0\right)  =0,\tag{4}\\
\tau_{4}u  &  :=u^{^{\prime}}\left(  h_{1}+0\right)  -\alpha_{4}u\left(
h_{1}-0\right)  -\beta_{4}u^{^{\prime}}\left(  h_{1}-0\right)  =0,\tag{5}\\
\tau_{5}u  &  :=u\left(  h_{2}+0\right)  -\alpha_{5}u\left(  h_{2}-0\right)
-\beta_{5}u^{^{\prime}}\left(  h_{2}-0\right)  =0,\tag{6}\\
\tau_{6}u  &  :=u^{^{\prime}}\left(  h_{2}+0\right)  -\alpha_{6}u\left(
h_{2}-0\right)  -\beta_{6}u^{^{\prime}}\left(  h_{2}-0\right)  =0,\tag{7}\\
\tau_{7}u  &  :=u\left(  h_{3}+0\right)  -\alpha_{7}u\left(  h_{3}-0\right)
-\beta_{7}u^{^{\prime}}\left(  h_{3}-0\right)  =0,\tag{8}\\
\tau_{8}u  &  :=u^{^{\prime}}\left(  h_{3}+0\right)  -\alpha_{8}u\left(
h_{3}-0\right)  -\beta_{8}u^{^{\prime}}\left(  h_{3}-0\right)  =0, \tag{9}%
\end{align}
where the coefficients $\alpha_{i},\beta_{i}$ and $\beta_{j}^{^{\prime}%
}\left(  i=\overline{1,8},\text{ }j=1,2\right)  $ are real numbers. Throughout
this paper, we assume that%

\[
\theta=\left\vert
\begin{array}
[c]{cc}%
\alpha_{3} & \beta_{3}\\
\alpha_{4} & \beta_{4}%
\end{array}
\right\vert >0,\gamma=\left\vert
\begin{array}
[c]{cc}%
\alpha_{5} & \beta_{5}\\
\alpha_{6} & \beta_{6}%
\end{array}
\right\vert >0,\xi=\left\vert
\begin{array}
[c]{cc}%
\alpha_{7} & \beta_{7}\\
\alpha_{8} & \beta_{8}%
\end{array}
\right\vert >0,
\]
\[
\rho=\left\vert
\begin{array}
[c]{cc}%
\beta_{1}^{^{\prime}} & \beta_{1}\\
\beta_{2}^{^{\prime}} & \beta_{2}%
\end{array}
\right\vert >0,
\]

and $\left\vert \alpha_{1}\right\vert +\left\vert \alpha_{2}\right\vert \neq0$.

\section{Operator formulation}

The relation between a symmetric linear operator $A$ defined in a suitable
Hilbert space $H$ and the problem $\left(  1\right)  -\left(  9\right)  $ has
been introduced in $\left[  13\right]  $. Here, we repeat the definition and
prove that the operator $A$ is self-adjoint, not only symmetric.

We define the inner product in $L^{2}\left(  I\right)  $ as%
\[
\left\langle f,g\right\rangle _{1}=\frac{1}{p_{1}^{2}}%
%TCIMACRO{\dint \limits_{-1}^{h_{1}}}%
%BeginExpansion
{\displaystyle\int\limits_{-1}^{h_{1}}}
%EndExpansion
f_{1}\overline{g_{1}}+\frac{1}{p_{2}^{2}\theta}%
%TCIMACRO{\dint \limits_{h_{1}}^{h_{2}}}%
%BeginExpansion
{\displaystyle\int\limits_{h_{1}}^{h_{2}}}
%EndExpansion
f_{2}\overline{g_{2}}+\frac{1}{p_{3}^{2}\theta\gamma}%
%TCIMACRO{\dint \limits_{h_{2}}^{1}}%
%BeginExpansion
{\displaystyle\int\limits_{h_{2}}^{1}}
%EndExpansion
f_{3}\overline{g_{3}}+\frac{1}{p_{4}^{2}\theta\gamma\xi}%
%TCIMACRO{\dint \limits_{h_{2}}^{1}}%
%BeginExpansion
{\displaystyle\int\limits_{h_{2}}^{1}}
%EndExpansion
f_{4}\overline{g_{4}},\text{\ \ }\forall f,g\in L^{2}\left(  I\right)  ,
\]
where%

\begin{align*}
f_{1}\left(  x\right)   &  :=\left\{
\begin{array}
[c]{ll}%
f(x), & x\in\lbrack-1,{h}_{1}),\\
\lim_{x\rightarrow h_{1}-0}f(x), & x=h_{1},
\end{array}
\right.  \text{ \ }f_{2}\left(  x\right)  :=\left\{
\begin{array}
[c]{cc}%
\lim_{x\rightarrow h_{1}+0}f(x), & x=h_{1},\\
f(x), & x\in\left(  h_{1},{h}_{2}\right)  ,\\
\lim_{x\rightarrow h_{2}-0}f(x), & x=h_{2},
\end{array}
\right. \\
f_{3}\left(  x\right)   &  :=\left\{
\begin{array}
[c]{cc}%
\lim_{x\rightarrow h_{2}+0}f(x), & x=h_{2},\\
f(x), & x\in\left(  h_{2},{h}_{3}\right)  ,\\
\lim_{x\rightarrow h_{3}-0}f(x), & x=h_{3},
\end{array}
\right.  \text{ \ }f_{4}\left(  x\right)  :=\left\{
\begin{array}
[c]{cc}%
\lim_{x\rightarrow h_{3}+0}f(x), & x=h_{3},\\
f\left(  x\right)  , & x\in\left(  h_{3},1\right]
\end{array}
\right.
\end{align*}
It is easy to verify that $\left(  L^{2}\left(  I\right)  ,\left\langle
\cdot,\cdot\right\rangle _{1}\right)  $ is a Hilbert space. For simplicity, it
is denoted by $H_{1}$.

The inner product in $H:=H_{1}\oplus%
%TCIMACRO{\U{2102} }%
%BeginExpansion
\mathbb{C}
%EndExpansion
$ is defined by%
\[
\left\langle F,G\right\rangle =\left\langle f,g\right\rangle _{1}+\frac
{1}{\rho\theta\gamma\xi}h\overline{k}%
\]
for $F=\left(  f\left(  x\right)  ,h\right)  $, $G=\left(  g\left(  x\right)
,k\right)  \in H$, where $f,g\in H_{1}$, $h,k\in%
%TCIMACRO{\U{2102} }%
%BeginExpansion
\mathbb{C}
%EndExpansion
$.

We define the operator $A$ in $H$ as follows:%
\[%
\begin{array}
[c]{c}%
D\left(  A\right)  :=\left\{  \left(  f\left(  x\right)  ,h\right)  \in H\mid
f_{1},f_{1}^{^{\prime}}\in AC_{loc}\left(  \left(  -1,h_{1}\right)  \right)
,f_{2},f_{2}^{^{\prime}}\in AC_{loc}\left(  \left(  h_{1},h_{2}\right)
\right)  ,\right. \\
f_{3},f_{3}^{^{\prime}}\in AC_{loc}\left(  \left(  h_{2},h_{3}\right)
\right)  ,\text{ }f_{4},f_{4}^{^{\prime}}\in AC_{loc}\left(  \left(
h_{3},1\right)  \right)  ,\text{ }\tau f\in H_{1},\tau_{1}f=\tau_{3}f=\tau
_{4}f=\tau_{5}f=\tau_{6}f,\\
\left.  h=\beta_{1}^{^{\prime}}f\left(  1\right)  -\beta_{2}^{^{\prime}%
}f^{^{\prime}}\left(  1\right)  ,\right\}
\end{array}
\]
$AF=\left(  \tau f,-\left(  \beta_{1}f\left(  1\right)  -\beta_{2}f^{^{\prime
}}\left(  1\right)  \right)  \right)  $ for\ $F=\left(  f,\beta_{1}^{^{\prime
}}f\left(  1\right)  -\beta_{2}^{^{\prime}}f\left(  1\right)  \right)  \in
D\left(  A\right)  $.

Note that by our assumption on $q\left(  x\right)  $ and Theorem 3.2 in
$\left[  19\right]  $, for each $\left(  f,h\right)  \in D\left(  A\right)  $,
$f_{1}$, $f_{1}^{^{\prime}}$ are continuous on $\left[  -1,h_{1}\right]  $,
$f_{2}$, $f_{2}^{^{\prime}}$ are continuous on $\left[  h_{1},h_{2}\right]  ,$
$f_{3}$, $f_{3}^{^{\prime}}$ are continuous on $\left[  h_{2},h_{3}\right]  $
and $f_{4}$, $f_{4}^{^{\prime}}$ are continuous on $\left[  h_{3},1\right]  $.
For simplicity, for $\left(  f,h\right)  \in D\left(  A\right)  ,$ set%
\[
N\left(  f\right)  =\beta_{1}f\left(  1\right)  -\beta_{2}f^{^{\prime}}\left(
1\right)  \text{, }N^{^{\prime}}\left(  f\right)  =\beta_{1}^{^{\prime}%
}f\left(  1\right)  -\beta_{2}^{^{\prime}}f^{^{\prime}}\left(  1\right)
\text{.}%
\]
So, we can study the problem $\left(  1\right)  -\left(  9\right)  $ in $H$ by
considering the operator equation $AF=\lambda F$. Obviously, we have

\textbf{Lemma 1.1 }\textit{The eigenvalues of the boundary value problem
}$\left(  1\right)  -\left(  9\right)  $ coincide with those of $A$\textit{,
and its eigenfunctions are the first components of the corresponding
eigenfunctions of }$A$.

\textbf{Lemma 1.2 }\textit{The domain }$D\left(  A\right)  $ is dense in $H$.

\begin{proof}
Suppose that $F\in H$ is orthogonal to all $G\in D\left(  A\right)  $ with
respect to the inner product $\left\langle \cdot,\cdot\right\rangle $, where
$F=\left(  f\left(  x\right)  ,h\right)  $, $G=\left(  g\left(  x\right)
,k\right)  $. Let $\widetilde{C}_{0}^{\infty}$ denote the set of functions%
\[
\phi\left(  x\right)  =\left\{
\begin{array}
[c]{c}%
\varphi_{1}\left(  x\right)  ,\text{ }x\in\left[  -1,h_{1}\right)  ,\\
\varphi_{2}\left(  x\right)  ,\text{ }x\in\left(  h_{1},{h}_{2}\right)  ,\\
\varphi_{3}\left(  x\right)  ,\text{ }x\in\left(  h_{2},{h}_{3}\right)  ,\\
\varphi_{4}\left(  x\right)  ,\text{ }x\in\left(  h_{3},{1}\right]  ,
\end{array}
\right.
\]
where $\varphi_{1}\left(  x\right)  \in C_{0}^{\infty}\left[  -1,h_{1}\right)
$, $\varphi_{2}\left(  x\right)  \in C_{0}^{\infty}\left(  h_{1},{h}%
_{2}\right)  ,$ $\varphi_{3}\left(  x\right)  \in C_{0}^{\infty}\left(
h_{2},{h}_{3}\right)  $ and $\varphi_{4}\left(  x\right)  \in C_{0}^{\infty
}\left(  h_{3},{1}\right]  $. Since $\widetilde{C}_{0}^{\infty}\oplus0\subset
D\left(  A\right)  $ $\left(  0\in%
%TCIMACRO{\U{2102} }%
%BeginExpansion
\mathbb{C}
%EndExpansion
\right)  $, any $U=\left(  u\left(  x\right)  ,0\right)  \in\widetilde{C}%
_{0}^{\infty}\oplus0$ is orthogonal to $F$, namely,%
\[
\left\langle F,U\right\rangle =\frac{1}{p_{1}^{2}}%
%TCIMACRO{\dint \limits_{-1}^{h_{1}}}%
%BeginExpansion
{\displaystyle\int\limits_{-1}^{h_{1}}}
%EndExpansion
f\left(  x\right)  \overline{u\left(  x\right)  }dx+\frac{1}{p_{2}^{2}\theta}%
%TCIMACRO{\dint \limits_{h_{1}}^{h_{2}}}%
%BeginExpansion
{\displaystyle\int\limits_{h_{1}}^{h_{2}}}
%EndExpansion
f\left(  x\right)  \overline{u\left(  x\right)  }dx+\frac{1}{p_{3}^{2}%
\theta\gamma}%
%TCIMACRO{\dint \limits_{h_{2}}^{h_{3}}}%
%BeginExpansion
{\displaystyle\int\limits_{h_{2}}^{h_{3}}}
%EndExpansion
f\left(  x\right)  \overline{u\left(  x\right)  }dx
\]%
\[
+\frac{1}{p_{4}^{2}\theta\gamma\xi}%
%TCIMACRO{\dint \limits_{h_{3}}^{1}}%
%BeginExpansion
{\displaystyle\int\limits_{h_{3}}^{1}}
%EndExpansion
f\left(  x\right)  \overline{u\left(  x\right)  }dx=\left\langle
f,u\right\rangle _{1}\text{.}%
\]
This implies that $f\left(  x\right)  $ is orthogonal to $\widetilde{C}%
_{0}^{\infty}$ in $H_{1}$ and hence vanishes. So, $\left\langle
F,G\right\rangle =\frac{1}{\rho\theta\gamma\xi}h\overline{k}=0$. Thus, $h=0$
since $k=N^{^{\prime}}\left(  g\right)  $ can be chosen arbitrarily. So,
$F=\left(  0,0\right)  $. Therefore, $D\left(  A\right)  $ is dense in $H.$
\end{proof}

\textbf{Theorem 1.1.} \textit{The linear operator }$A$\textit{ is self-adjoint
in }$H$\textit{.}

\begin{proof}
For all $F,G\in D\left(  A\right)  $, (2) implies that $f\left(  -1\right)
\overline{g}^{^{\prime}}\left(  -1\right)  -f^{^{\prime}}\left(  -1\right)
\overline{g}\left(  -1\right)  =0$, and direct calculations using (4) and (5)
then yield that%
\begin{align*}
\left\langle AF,G\right\rangle  &  =\left\langle F,AG\right\rangle +W\left(
f,\overline{g};h_{1}-0\right)  -W\left(  f,\overline{g};-1\right)  +\frac
{1}{\theta}W\left(  f,\overline{g};h_{2}-0\right)  -\\
&  \frac{1}{\theta}W\left(  f,\overline{g};h_{1}+0\right)  +\frac{1}%
{\theta\gamma}W\left(  f,\overline{g};h_{3}-0\right)  -\frac{1}{\theta\gamma
}W\left(  f,\overline{g};h_{2}+0\right)  +\\
&  \frac{1}{\theta\gamma\xi}W\left(  f,\overline{g};1\right)  -\frac{1}%
{\theta\gamma\xi}W\left(  f,\overline{g};h_{3}+0\right)  -\frac{1}{\rho
\theta\gamma\xi}\left(  N\left(  f\right)  \overline{N^{^{\prime}}\left(
g\right)  }-N^{^{\prime}}\left(  f\right)  \overline{N\left(  g\right)
}\right) \\
&  =\left\langle F,AG\right\rangle ,
\end{align*}
where $\left(  f,g;x\right)  $ denotes the Wronskians $f\left(  x\right)
g^{^{\prime}}\left(  x\right)  -f^{^{\prime}}\left(  x\right)  g\left(
x\right)  $. So, $A$ is symmetric.

It remains to show that if $\left\langle AF,W\right\rangle =\left\langle
F,U\right\rangle $ for all $F=\left(  f,N^{\prime}\left(  f\right)  \right)
\in D\left(  A\right)  $, then $W\in D\left(  A\right)  $ and $AW=U$, where
$W=\left(  w\left(  x\right)  ,h\right)  $ and $U=\left(  u\left(  x\right)
,k\right)  $, i.e., $\left(  i\right)  $ $w_{1},w_{1}^{^{\prime}}\in
AC_{loc}\left(  \left(  -1,h_{1}\right)  \right)  $, $w_{2},w_{2}^{^{\prime}%
}\in AC_{loc}\left(  \left(  h_{1},h_{2}\right)  \right)  $, $w_{3}%
,w_{3}^{^{\prime}}\in AC_{loc}\left(  \left(  h_{2},h_{3}\right)  \right)  $
$w_{4},w_{4}^{^{\prime}}\in AC_{loc}\left(  \left(  h_{3},1\right)  \right)  $
and $\tau w\in H_{1};$ $\left(  ii\right)  $ $h=N^{^{\prime}}\left(  w\right)
=\beta_{1}^{^{\prime}}w\left(  1\right)  -\beta_{2}^{^{\prime}}w^{^{\prime}%
}\left(  1\right)  ;$ $\left(  iii\right)  $ $\tau_{1}w=\tau_{3}w=\tau
_{4}w=\tau_{5}w=\tau_{6}w=\tau_{7}w=\tau_{8}w=0;$ $\left(  iv\right)  $
$u\left(  x\right)  =\tau w;$ $\left(  v\right)  $ $k=-N\left(  w\right)
=-\beta_{1}w\left(  1\right)  +\beta_{2}w^{^{\prime}}\left(  1\right)  $.

For all $F\in\widetilde{C}_{0}^{\infty}\oplus0\subset D\left(  A\right)  $, we
obtain%
\begin{align*}
&  \frac{1}{p_{1}^{2}}%
%TCIMACRO{\dint \limits_{-1}^{h_{1}}}%
%BeginExpansion
{\displaystyle\int\limits_{-1}^{h_{1}}}
%EndExpansion
\left(  \tau f\right)  \overline{w}dx+\frac{1}{p_{2}^{2}\theta}%
%TCIMACRO{\dint \limits_{h_{1}}^{h_{2}}}%
%BeginExpansion
{\displaystyle\int\limits_{h_{1}}^{h_{2}}}
%EndExpansion
\left(  \tau f\right)  \overline{w}dx+\frac{1}{p_{3}^{2}\theta\gamma}%
%TCIMACRO{\dint \limits_{h_{2}}^{h_{3}}}%
%BeginExpansion
{\displaystyle\int\limits_{h_{2}}^{h_{3}}}
%EndExpansion
\left(  \tau f\right)  \overline{w}dx+\frac{1}{p_{4}^{2}\theta\gamma\xi}%
%TCIMACRO{\dint \limits_{h_{3}}^{1}}%
%BeginExpansion
{\displaystyle\int\limits_{h_{3}}^{1}}
%EndExpansion
\left(  \tau f\right)  \overline{w}dx\\
&  =\frac{1}{p_{1}^{2}}%
%TCIMACRO{\dint \limits_{-1}^{h_{1}}}%
%BeginExpansion
{\displaystyle\int\limits_{-1}^{h_{1}}}
%EndExpansion
f\overline{u}dx+\frac{1}{p_{2}^{2}\theta}%
%TCIMACRO{\dint \limits_{h_{1}}^{h_{2}}}%
%BeginExpansion
{\displaystyle\int\limits_{h_{1}}^{h_{2}}}
%EndExpansion
f\overline{u}dx+\frac{1}{p_{3}^{2}\theta\gamma}%
%TCIMACRO{\dint \limits_{h_{2}}^{h_{3}}}%
%BeginExpansion
{\displaystyle\int\limits_{h_{2}}^{h_{3}}}
%EndExpansion
f\overline{u}dx+\frac{1}{p_{4}^{2}\theta\gamma\xi}%
%TCIMACRO{\dint \limits_{h_{3}}^{1}}%
%BeginExpansion
{\displaystyle\int\limits_{h_{3}}^{1}}
%EndExpansion
f\overline{u}dx\text{,}%
\end{align*}
namely, $\left\langle \tau f,w\right\rangle _{1}=\left\langle f,u\right\rangle
_{1}$. Hence, by standart Sturm-Liouville theory, $\left(  i\right)  $ and
$\left(  iv\right)  $ hold. By $\left(  iv\right)  $, the equation
$\left\langle AF,W\right\rangle =\left\langle F,U\right\rangle $, $\forall
F\in D\left(  A\right)  $, becomes%
\begin{align*}
&  \frac{1}{p_{1}^{2}}%
%TCIMACRO{\dint \limits_{-1}^{h_{1}}}%
%BeginExpansion
{\displaystyle\int\limits_{-1}^{h_{1}}}
%EndExpansion
\left(  \tau f\right)  \overline{w}dx+\frac{1}{p_{2}^{2}\theta}%
%TCIMACRO{\dint \limits_{h_{1}}^{h_{2}}}%
%BeginExpansion
{\displaystyle\int\limits_{h_{1}}^{h_{2}}}
%EndExpansion
\left(  \tau f\right)  \overline{w}dx+\frac{1}{p_{3}^{2}\theta\gamma}%
%TCIMACRO{\dint \limits_{h_{2}}^{h_{3}}}%
%BeginExpansion
{\displaystyle\int\limits_{h_{2}}^{h_{3}}}
%EndExpansion
\left(  \tau f\right)  \overline{w}dx+\frac{1}{p_{4}^{2}\theta\gamma\xi}%
%TCIMACRO{\dint \limits_{h_{3}}^{1}}%
%BeginExpansion
{\displaystyle\int\limits_{h_{3}}^{1}}
%EndExpansion
\left(  \tau f\right)  \overline{w}dx-\frac{N\left(  f\right)  \overline{h}%
}{\rho\theta\gamma}\\
&  =\frac{1}{p_{1}^{2}}%
%TCIMACRO{\dint \limits_{-1}^{h_{1}}}%
%BeginExpansion
{\displaystyle\int\limits_{-1}^{h_{1}}}
%EndExpansion
f\left(  \tau\overline{w}\right)  dx+\frac{1}{p_{2}^{2}\theta}%
%TCIMACRO{\dint \limits_{h_{1}}^{h_{2}}}%
%BeginExpansion
{\displaystyle\int\limits_{h_{1}}^{h_{2}}}
%EndExpansion
f\left(  \tau\overline{w}\right)  dx+\frac{1}{p_{3}^{2}\theta\gamma}%
%TCIMACRO{\dint \limits_{h_{2}}^{h_{3}}}%
%BeginExpansion
{\displaystyle\int\limits_{h_{2}}^{h_{3}}}
%EndExpansion
f\left(  \tau\overline{w}\right)  dx+\frac{1}{p_{4}^{2}\theta\gamma\xi}%
%TCIMACRO{\dint \limits_{h_{3}}^{1}}%
%BeginExpansion
{\displaystyle\int\limits_{h_{3}}^{1}}
%EndExpansion
f\left(  \tau\overline{w}\right)  dx+\frac{N^{^{\prime}}\left(  f\right)
\overline{k}}{\rho\theta\gamma\xi}.
\end{align*}
So,%
\[
\left\langle \tau f,w\right\rangle _{1}=\left\langle f,\tau w\right\rangle
_{1}+\frac{N^{^{\prime}}\left(  f\right)  \overline{k}+N\left(  f\right)
\overline{h}}{\rho\theta\gamma\xi}.
\]
However,%
\begin{align*}
\left\langle \tau f,w\right\rangle _{1}  &  =\frac{1}{p_{1}^{2}}%
%TCIMACRO{\dint \limits_{-1}^{h_{1}}}%
%BeginExpansion
{\displaystyle\int\limits_{-1}^{h_{1}}}
%EndExpansion
\left(  -p_{1}^{2}f^{^{\prime\prime}}+q\left(  x\right)  f\right)
\overline{w}dx+\frac{1}{p_{2}^{2}\theta}%
%TCIMACRO{\dint \limits_{h_{1}}^{h_{2}}}%
%BeginExpansion
{\displaystyle\int\limits_{h_{1}}^{h_{2}}}
%EndExpansion
\left(  -p_{2}^{2}f^{^{\prime\prime}}+q\left(  x\right)  f\right)
\overline{w}dx+\\
&  \frac{1}{p_{3}^{2}\theta\gamma}%
%TCIMACRO{\dint \limits_{h_{2}}^{h_{3}}}%
%BeginExpansion
{\displaystyle\int\limits_{h_{2}}^{h_{3}}}
%EndExpansion
\left(  -p_{3}^{2}f^{^{\prime\prime}}+q\left(  x\right)  f\right)
\overline{w}dx+\frac{1}{p_{4}^{2}\theta\gamma\xi}%
%TCIMACRO{\dint \limits_{h_{3}}^{1}}%
%BeginExpansion
{\displaystyle\int\limits_{h_{3}}^{1}}
%EndExpansion
\left(  -p_{4}^{2}f^{^{\prime\prime}}+q\left(  x\right)  f\right)
\overline{w}dx\\
&  =\frac{1}{p_{1}^{2}}%
%TCIMACRO{\dint \limits_{-1}^{h_{1}}}%
%BeginExpansion
{\displaystyle\int\limits_{-1}^{h_{1}}}
%EndExpansion
f\left(  \tau\overline{w}\right)  dx+\frac{1}{p_{2}^{2}\theta}%
%TCIMACRO{\dint \limits_{h_{1}}^{h_{2}}}%
%BeginExpansion
{\displaystyle\int\limits_{h_{1}}^{h_{2}}}
%EndExpansion
f\left(  \tau\overline{w}\right)  dx+\frac{1}{p_{3}^{2}\theta\gamma}%
%TCIMACRO{\dint \limits_{h_{2}}^{h_{3}}}%
%BeginExpansion
{\displaystyle\int\limits_{h_{2}}^{h_{3}}}
%EndExpansion
f\left(  \tau\overline{w}\right)  dx+\frac{1}{p_{4}^{2}\theta\gamma\xi}%
%TCIMACRO{\dint \limits_{h_{3}}^{1}}%
%BeginExpansion
{\displaystyle\int\limits_{h_{3}}^{1}}
%EndExpansion
f\left(  \tau\overline{w}\right)  dx\\
+  &  W\left(  f,\overline{w};h_{1}-0\right)  -W\left(  f,\overline
{w};-1\right)  +\frac{1}{\theta}W\left(  f,\overline{w};h_{2}-0\right)
-\frac{1}{\theta}W\left(  f,\overline{w};h_{1}+0\right)  +\\
&  \frac{1}{\theta\gamma}W\left(  f,\overline{w};h_{3}-0\right)  -\frac
{1}{\theta\gamma}W\left(  f,\overline{w};h_{2}+0\right)  +\frac{1}%
{\theta\gamma\xi}W\left(  f,\overline{w};1\right)  -\frac{1}{\theta\gamma\xi
}W\left(  f,\overline{w};h_{3}+0\right) \\
&  =\left\langle f,\tau w\right\rangle _{1}+W\left(  f,\overline{w}%
;h_{1}-0\right)  -W\left(  f,\overline{w};-1\right)  +\frac{1}{\theta}W\left(
f,\overline{w};h_{2}-0\right)  -\frac{1}{\theta}W\left(  f,\overline{w}%
;h_{1}+0\right) \\
&  +\frac{1}{\theta\gamma}W\left(  f,\overline{w};h_{3}-0\right)  -\frac
{1}{\theta\gamma}W\left(  f,\overline{w};h_{2}+0\right)  +\frac{1}%
{\theta\gamma\xi}W\left(  f,\overline{w};1\right)  -\frac{1}{\theta\gamma\xi
}W\left(  f,\overline{w};h_{3}+0\right)  .
\end{align*}
Hence,%
\begin{align}
\frac{N^{^{\prime}}\left(  f\right)  \overline{k}+N\left(  f\right)
\overline{h}}{\rho\theta\gamma\xi}  &  =W\left(  f,\overline{w};h_{1}%
-0\right)  -W\left(  f,\overline{w};-1\right)  +\frac{1}{\theta}W\left(
f,\overline{w};h_{2}-0\right) \nonumber\\
&  -\frac{1}{\theta}W\left(  f,\overline{w};h_{1}+0\right)  +\frac{1}%
{\theta\gamma}W\left(  f,\overline{w};h_{3}-0\right)  -\frac{1}{\theta\gamma
}W\left(  f,\overline{w};h_{2}+0\right) \nonumber\\
&  +\frac{1}{\theta\gamma\xi}W\left(  f,\overline{w};1\right)  -\frac
{1}{\theta\gamma\xi}W\left(  f,\overline{w};h_{2}+0\right) \nonumber\\
&  =\left(  f\left(  h_{1}-0\right)  \overline{w}^{^{\prime}}\left(
h_{1}-0\right)  -f^{^{\prime}}\left(  h_{1}-0\right)  \overline{w}\left(
h_{1}-0\right)  \right)  -\left(  f\left(  -1\right)  \overline{w}^{^{\prime}%
}\left(  -1\right)  \right. \nonumber\\
&  \left.  -f^{^{\prime}}\left(  -1\right)  \overline{w}\left(  -1\right)
\right)  +\frac{1}{\theta}\left(  f\left(  h_{2}-0\right)  \overline
{w}^{^{\prime}}\left(  h_{2}-0\right)  -f^{^{\prime}}\left(  h_{2}-0\right)
\overline{w}\left(  h_{2}-0\right)  \right) \nonumber\\
&  -\frac{1}{\theta}\left(  f\left(  h_{1}+0\right)  \overline{w}^{^{\prime}%
}\left(  h_{1}+0\right)  -f^{^{\prime}}\left(  h_{1}+0\right)  \overline
{w}\left(  h_{1}+0\right)  \right) \nonumber\\
&  +\frac{1}{\theta\gamma}\left(  f\left(  h_{3}-0\right)  \overline
{w}^{^{\prime}}\left(  h_{3}-0\right)  -f^{^{\prime}}\left(  h_{3}-0\right)
\overline{w}\left(  h_{3}-0\right)  \right) \nonumber\\
&  -\frac{1}{\theta\gamma}\left(  f\left(  h_{2}+0\right)  \overline
{w}^{^{\prime}}\left(  h_{2}+0\right)  -f^{^{\prime}}\left(  h_{2}+0\right)
\overline{w}\left(  h_{2}+0\right)  \right)  +\frac{1}{\theta\gamma\xi}\left(
f\left(  1\right)  \overline{w}^{^{\prime}}\left(  1\right)  \right.
\nonumber\\
&  \left.  -f^{^{\prime}}\left(  1\right)  \overline{w}\left(  1\right)
\right)  -\frac{1}{\theta\gamma\xi}\left(  f\left(  h_{3}+0\right)
\overline{w}^{^{\prime}}\left(  h_{3}+0\right)  -f^{^{\prime}}\left(
h_{3}+0\right)  \overline{w}\left(  h_{3}+0\right)  \right)  . \tag{8}%
\end{align}
By Naimark's Patching Lemma $\left[  20\right]  $, there is an $F\in D\left(
A\right)  $ such that $f\left(  -1\right)  =f^{^{\prime}}\left(  -1\right)
=f\left(  h_{1}-0\right)  =f^{^{\prime}}\left(  h_{1}-0\right)  =f\left(
h_{1}+0\right)  =f^{^{\prime}}\left(  h_{1}+0\right)  =f\left(  h_{2}%
-0\right)  =f^{^{\prime}}\left(  h_{2}-0\right)  =f\left(  h_{2}+0\right)
=f^{^{\prime}}\left(  h_{2}+0\right)  =f\left(  h_{3}-0\right)  =f^{^{\prime}%
}\left(  h_{3}-0\right)  =f\left(  h_{3}+0\right)  =f^{^{\prime}}\left(
h_{3}+0\right)  =0$, $f\left(  1\right)  =\beta_{2}^{^{\prime}}$ and
$f^{^{\prime}}\left(  1\right)  =\beta_{1}^{^{\prime}}$. For such an $F$,
$N^{^{\prime}}\left(  f\right)  =0$. Thus, from $\left(  8\right)  $ we obtain
$h=\beta_{1}^{^{\prime}}w\left(  1\right)  -\beta_{2}^{^{\prime}}w^{^{\prime}%
}\left(  1\right)  $. Namely, $\left(  ii\right)  $ holds. Similarly, one
proves $\left(  v\right)  $.

It remains to show that $\left(  iii\right)  $ holds. Choose $F\in D\left(
A\right)  $ so that $f\left(  1\right)  =f^{^{\prime}}\left(  1\right)
=f\left(  h_{1}-0\right)  =f^{^{\prime}}\left(  h_{1}-0\right)  =f\left(
h_{2}-0\right)  =f^{^{\prime}}\left(  h_{2}-0\right)  =f\left(  h_{3}%
-0\right)  =f^{^{\prime}}\left(  h_{3}-0\right)  =0,$ $f\left(  -1\right)
=\alpha_{2}$ and $f^{^{\prime}}\left(  -1\right)  =-\alpha_{1}$. $N^{^{\prime
}}\left(  f\right)  =N\left(  f\right)  =0.$ From $\left(  8\right)  ,$ we get
$\alpha_{1}w\left(  -1\right)  +\alpha_{2}w^{^{\prime}}\left(  -1\right)  =0$.
Let $F\in D\left(  A\right)  $ satisfies $f\left(  1\right)  =f^{^{\prime}%
}\left(  1\right)  =f\left(  -1\right)  =f^{^{\prime}}\left(  -1\right)
=f\left(  h_{1}+0\right)  =f\left(  h_{2}+0\right)  =f\left(  h_{3}+0\right)
=0,$ $f\left(  h_{1}-0\right)  =-\beta_{3},$ $f\left(  h_{2}-0\right)
=-\beta_{5},$ $f\left(  h_{3}-0\right)  =-\beta_{7},$ $f^{^{\prime}}\left(
h_{1}-0\right)  =\alpha_{3},$ $f^{^{\prime}}\left(  h_{2}-0\right)
=\alpha_{5},$ $f^{^{\prime}}\left(  h_{3}-0\right)  =\alpha_{7},$
$f^{^{\prime}}\left(  h_{1}+0\right)  =\theta,$ $f^{^{\prime}}\left(
h_{2}+0\right)  =\gamma$ and $f^{^{\prime}}\left(  h_{3}+0\right)  =\xi$. Then
$N\left(  f\right)  =N^{^{\prime}}\left(  f\right)  =0.$ By $\left(  8\right)
$, we have $w\left(  h_{1}+0\right)  =\alpha_{3}w\left(  h_{1}-0\right)
+\beta_{3}w^{^{\prime}}\left(  h_{1}-0\right)  ,$ $w\left(  h_{2}+0\right)
=\alpha_{5}w\left(  h_{2}-0\right)  +\beta_{5}w^{^{\prime}}\left(
h_{2}-0\right)  $ and $w\left(  h_{3}+0\right)  =\alpha_{7}w\left(
h_{3}-0\right)  +\beta_{7}w^{^{\prime}}\left(  h_{3}-0\right)  .$ Finally,
choose $F\in D\left(  A\right)  $ so that $f\left(  1\right)  =f^{^{\prime}%
}\left(  1\right)  =f\left(  -1\right)  =f^{^{\prime}}\left(  -1\right)
=f^{^{\prime}}\left(  h_{1}+0\right)  =f^{^{\prime}}\left(  h_{2}+0\right)
=f^{^{\prime}}\left(  h_{3}+0\right)  =0,$ $f\left(  h_{1}-0\right)
=\beta_{4},$ $f^{^{\prime}}\left(  h_{1}-0\right)  =-\alpha_{4},$ $f\left(
h_{2}-0\right)  =\beta_{6},$ $f^{^{\prime}}\left(  h_{2}-0\right)
=-\alpha_{6},$ $f\left(  h_{3}-0\right)  =\beta_{8},$ $f^{^{\prime}}\left(
h_{3}-0\right)  =-\alpha_{8}$ $f\left(  h_{1}+0\right)  =\theta,$ $f\left(
h_{2}+0\right)  =\gamma$ and $f\left(  h_{3}+0\right)  =\xi.$ Then $N\left(
f\right)  =N^{^{\prime}}\left(  f\right)  =0.$ From $\left(  8\right)  ,$ we
obtain $w^{^{\prime}}\left(  h_{1}+0\right)  =\alpha_{4}w\left(
h_{1}-0\right)  +\beta_{4}w^{^{\prime}}\left(  h_{1}-0\right)  $,
$w^{^{\prime}}\left(  h_{2}+0\right)  =\alpha_{6}w\left(  h_{2}-0\right)
+\beta_{6}w^{^{\prime}}\left(  h_{2}-0\right)  $ and $w^{^{\prime}}\left(
h_{3}+0\right)  =\alpha_{8}w\left(  h_{3}-0\right)  +\beta_{8}w^{^{\prime}%
}\left(  h_{3}-0\right)  .$
\end{proof}

\textbf{Corollary 1.1} \textit{The eigenvalues of }$\left(  1\right)  -\left(
9\right)  $\textit{ are real, and if }$\lambda_{1}$\textit{ and }$\lambda_{2}%
$\textit{ are two different eigenvalues of }$\left(  1\right)  -\left(
9\right)  $\textit{, then the corresponding eigenfunctions }$f\left(
x\right)  $\textit{ and }$g\left(  x\right)  $\textit{ are orthogonal in the
sense of}%
\[
\frac{1}{p_{1}^{2}}%
%TCIMACRO{\dint \limits_{-1}^{h_{1}}}%
%BeginExpansion
{\displaystyle\int\limits_{-1}^{h_{1}}}
%EndExpansion
f\overline{g}+\frac{1}{p_{2}^{2}\theta}%
%TCIMACRO{\dint \limits_{h_{1}}^{h_{2}}}%
%BeginExpansion
{\displaystyle\int\limits_{h_{1}}^{h_{2}}}
%EndExpansion
f\overline{g}+\frac{1}{p_{3}^{2}\theta\gamma}%
%TCIMACRO{\dint \limits_{h_{2}}^{h_{3}}}%
%BeginExpansion
{\displaystyle\int\limits_{h_{2}}^{h_{3}}}
%EndExpansion
f\overline{g}+\frac{1}{p_{4}^{2}\theta\gamma\xi}%
%TCIMACRO{\dint \limits_{h_{3}}^{1}}%
%BeginExpansion
{\displaystyle\int\limits_{h_{3}}^{1}}
%EndExpansion
f\overline{g}+
\]%
\[
\frac{1}{\rho\theta\gamma\xi}\left(  \beta_{1}^{^{\prime}}f\left(  1\right)
-\beta_{2}^{^{\prime}}f^{^{\prime}}\left(  1\right)  \right)  \left(
\beta_{1}^{^{\prime}}\overline{g}\left(  1\right)  -\beta_{2}^{^{\prime}%
}\overline{g}^{^{\prime}}\left(  1\right)  \right)  =0.
\]

\section{Simplicity of eigenvalues}

We consider the initial-value problem%
\[
\left\{
\begin{array}
[c]{c}%
-p_{1}^{2}u^{^{\prime\prime}}\left(  x\right)  +q(x)u\left(  x\right)
=\lambda u\left(  x\right)  ,\text{ \ \ }x\in\left[  -1,h_{1}\right)  ,\\
u\left(  -1\right)  =\alpha_{2},\text{ }u^{^{\prime}}\left(  -1\right)
=-\alpha_{1}.
\end{array}
\right.
\]
In terms of existence and uniqueness in ordinary differential equation theory,
the initial-value problem has a unique solution $\varphi_{1}\left(
x,\lambda\right)  $ for every $\lambda\in%
%TCIMACRO{\U{2102} }%
%BeginExpansion
\mathbb{C}
%EndExpansion
.$ Similarly, the initial-value problem%
\[
\left\{
\begin{array}
[c]{c}%
-p_{2}^{2}u^{^{\prime\prime}}\left(  x\right)  +q(x)u\left(  x\right)
=\lambda u\left(  x\right)  ,\text{ \ \ }x\in\left(  h_{1},h_{2}\right)  ,\\
u\left(  h_{1}\right)  =\alpha_{3}\varphi_{1}\left(  h_{1},\lambda\right)
+\beta_{3}\varphi_{1}^{^{\prime}}\left(  h_{1},\lambda\right)  ,\\
u^{^{\prime}}\left(  h_{1}\right)  =\alpha_{4}\varphi_{1}\left(  h_{1}%
,\lambda\right)  +\beta_{4}\varphi_{1}^{^{\prime}}\left(  h_{1},\lambda
\right)
\end{array}
\right.
\]
has a unique solution $\varphi_{2}\left(  x,\lambda\right)  $ for every
$\lambda\in%
%TCIMACRO{\U{2102} }%
%BeginExpansion
\mathbb{C}
%EndExpansion
.$The initial-value problem%
\[
\left\{
\begin{array}
[c]{c}%
-p_{3}^{2}u^{^{\prime\prime}}\left(  x\right)  +q(x)u\left(  x\right)
=\lambda u\left(  x\right)  ,\text{ \ \ }x\in\left(  h_{2},h_{3}\right)  ,\\
u\left(  h_{2}\right)  =\alpha_{5}\varphi_{2}\left(  h_{2},\lambda\right)
+\beta_{5}\varphi_{2}^{^{\prime}}\left(  h_{2},\lambda\right)  ,\\
u^{^{\prime}}\left(  h_{2}\right)  =\alpha_{6}\varphi_{2}\left(  h_{2}%
,\lambda\right)  +\beta_{6}\varphi_{2}^{^{\prime}}\left(  h_{2},\lambda
\right)
\end{array}
\right.
\]
has a unique solution $\varphi_{3}\left(  x,\lambda\right)  $ for every
$\lambda\in%
%TCIMACRO{\U{2102} }%
%BeginExpansion
\mathbb{C}
%EndExpansion
.$ Similarly, the initial-value problem%
\[
\left\{
\begin{array}
[c]{c}%
-p_{4}^{2}u^{^{\prime\prime}}\left(  x\right)  +q(x)u\left(  x\right)
=\lambda u\left(  x\right)  ,\text{ \ \ }x\in\left(  h_{3},1\right]  ,\\
u\left(  h_{3}\right)  =\alpha_{7}\varphi_{2}\left(  h_{2},\lambda\right)
+\beta_{5}\varphi_{2}^{^{\prime}}\left(  h_{2},\lambda\right)  ,\\
u^{^{\prime}}\left(  h_{3}\right)  =\alpha_{8}\varphi_{2}\left(  h_{2}%
,\lambda\right)  +\beta_{6}\varphi_{2}^{^{\prime}}\left(  h_{2},\lambda
\right)
\end{array}
\right.
\]
has a unique solution $\varphi_{4}\left(  x,\lambda\right)  $ for every
$\lambda\in%
%TCIMACRO{\U{2102} }%
%BeginExpansion
\mathbb{C}
%EndExpansion
$. For each given $x\in\left[  -1,h_{1}\right)  ,$ $\varphi_{1}\left(
x,\lambda\right)  $ is an entire function of $\lambda;$ for every $x\in\left(
h_{1},h_{2}\right)  ,$ $\varphi_{2}\left(  x,\lambda\right)  $ is an entire
function of $\lambda;$ for every $x\in\left(  h_{2},h_{3}\right)  ,$
$\varphi_{3}\left(  x,\lambda\right)  $ is an entire function of $\lambda$ and
for every $x\in\left(  h_{3},1\right]  ,$ $\varphi_{4}\left(  x,\lambda
\right)  $ is an entire function of $\lambda$.

Now we define a function $\phi\left(  x,\lambda\right)  $ on $x\in\left[
-1,h_{1}\right)  \cup\left(  h_{1},h_{2}\right)  \cup\left(  h_{2}%
,h_{3}\right)  \cup\left(  h_{3},1\right]  $ by%
\[
\phi\left(  x,\lambda\right)  =\left\{
\begin{array}
[c]{c}%
\varphi_{1}\left(  x,\lambda\right)  ,\text{ }x\in\left[  -1,h_{1}\right)  ,\\
\varphi_{2}\left(  x,\lambda\right)  ,\text{ }x\in\left(  h_{1},{h}%
_{2}\right)  ,\\
\varphi_{3}\left(  x,\lambda\right)  ,\text{ }x\in\left(  h_{2},{h}%
_{3}\right)  ,\\
\varphi_{4}\left(  x,\lambda\right)  ,\text{ }x\in\left(  h_{3},{1}\right]  .
\end{array}
\right.
\]
Obviously $\phi\left(  x,\lambda\right)  $ satisfies $\left(  1\right)  ,$
$\left(  2\right)  $ and $\left(  4\right)  -\left(  9\right)  .$ Similarly,
we define the function%
\[
\chi\left(  x,\lambda\right)  =\left\{
\begin{array}
[c]{c}%
\chi_{1}\left(  x,\lambda\right)  ,\text{ }x\in\left[  -1,h_{1}\right)  ,\\
\chi_{2}\left(  x,\lambda\right)  ,\text{ }x\in\left(  h_{1},{h}_{2}\right)
,\\
\chi_{3}\left(  x,\lambda\right)  ,\text{ }x\in\left(  h_{2},{h}_{3}\right)
,\\
\chi_{4}\left(  x,\lambda\right)  ,\text{ }x\in\left(  h_{3},{1}\right]  ,
\end{array}
\right.
\]
which satisfies $\left(  1\right)  ,$ $\left(  3\right)  -\left(  9\right)  .$

The Wronskian $W\left(  \varphi_{i}\left(  x,\lambda\right)  ,\chi_{i}\left(
x,\lambda\right)  \right)  $ $\left(  i=1,2,3,4\right)  $ are independent of
the variable $x$. Let $w_{i}\left(  \lambda\right)  =W\left(  \varphi
_{i}\left(  x,\lambda\right)  ,\chi_{i}\left(  x,\lambda\right)  \right)  $
and $w\left(  \lambda\right)  =w_{1}\left(  \lambda\right)  ,$ and then we
obtain $w_{2}\left(  \lambda\right)  =\theta w\left(  \lambda\right)  ,$
$w_{3}\left(  \lambda\right)  =\theta\gamma w\left(  \lambda\right)  $ and
$w_{4}\left(  \lambda\right)  =\theta\gamma\xi w\left(  \lambda\right)  .$

\textbf{Lemma 2.1 }$\left[  10\right]  $\textit{ The eigenvalues of the
problem }$\left(  1\right)  -\left(  9\right)  $\textit{ coincide with the
zeros of the entire function }$\mathit{w}\left(  \lambda\right)  $\textit{.}

\textbf{Definition 2.1 }\textit{The analytic multiplicity of an eigenvalue
}$\lambda$\textit{ of }$\left(  1\right)  -\left(  9\right)  $\textit{ is its
order as a root of the characteristic equation }$w\left(  \lambda\right)
=0.$\textit{ The geometric multiplicity of an eigenvalue is the dimension of
its eigenspace, i.e., the number of its linearly independent eigenfunctions.}

For convenience, set $\phi=\phi\left(  x,\lambda\right)  ,$ $\chi_{1\lambda
}=\frac{\partial\chi_{1}}{\partial\lambda},$ $\chi_{1\lambda}^{^{\prime}%
}=\frac{\partial\chi_{1}^{^{\prime}}}{\partial\lambda},$ etc.

\textbf{Theorem 2.1 }\textit{The eigenvalues of }$\left(  1\right)  -\left(
9\right)  $\textit{ are analytically simple.}

\begin{proof}
Let $\lambda=s+it$, where $s,t\in%
%TCIMACRO{\U{211d} }%
%BeginExpansion
\mathbb{R}
%EndExpansion
$ and $i^{2}=-1$. We differentiate the equation $\tau\chi=\lambda\chi$ with
respect to $\lambda$ and have%
\[
\tau\chi_{\lambda}=\lambda\chi_{\lambda}+\chi.
\]
By integration by parts, we get%
\begin{align}
\left\langle \tau\chi_{\lambda},\phi\right\rangle _{1}-\left\langle
\chi_{\lambda},\tau\phi\right\rangle _{1} &  =\left(  \chi_{1\lambda}%
\overline{\varphi_{1}}^{^{\prime}}-\chi_{1\lambda}^{^{\prime}}\overline
{\varphi_{1}}\right)  \mid_{-1}^{h_{1}}+\frac{1}{\theta}\left(  \chi
_{2\lambda}\overline{\varphi_{2}}^{^{\prime}}-\chi_{2\lambda}^{^{\prime}%
}\overline{\varphi_{2}}\right)  \mid_{h_{1}}^{h_{2}}+\nonumber\\
\frac{1}{\theta\gamma}\left(  \chi_{3\lambda}\overline{\varphi_{3}}^{^{\prime
}}-\chi_{3\lambda}^{^{\prime}}\overline{\varphi_{3}}\right)   &  \mid_{h_{2}%
}^{h_{3}}+\frac{1}{\theta\gamma\xi}\left(  \chi_{4\lambda}\overline
{\varphi_{4}}^{^{\prime}}-\chi_{4\lambda}^{^{\prime}}\overline{\varphi_{4}%
}\right)  \mid_{h_{3}}^{1}.\tag{9}%
\end{align}
Substituting $\tau\chi_{\lambda}=\lambda\chi_{\lambda}+\chi$ and $\tau
\phi=\lambda\phi$ into the left side of $\left(  9\right)  ,$ we have
\[
\lambda\left\langle \chi_{\lambda},\phi\right\rangle _{1}+\left\langle
\chi,\phi\right\rangle _{1}-\left\langle \chi_{\lambda},\lambda\phi
\right\rangle _{1}=\left\langle \chi,\phi\right\rangle _{1}+2it\left\langle
\chi_{\lambda},\phi\right\rangle _{1}.
\]
Moreover,%
\[
\left(  \chi_{1\lambda}\overline{\varphi_{1}}^{^{\prime}}-\chi_{1\lambda
}^{^{\prime}}\overline{\varphi_{1}}\right)  \mid_{-1}^{h_{1}}+\frac{1}{\theta
}\left(  \chi_{2\lambda}\overline{\varphi_{2}}^{^{\prime}}-\chi_{2\lambda
}^{^{\prime}}\overline{\varphi_{2}}\right)  \mid_{h_{1}}^{h_{2}}%
\]%
\[
+\frac{1}{\theta\gamma}\left(  \chi_{3\lambda}\overline{\varphi_{3}}%
^{^{\prime}}-\chi_{3\lambda}^{^{\prime}}\overline{\varphi_{3}}\right)
\mid_{h_{2}}^{h_{3}}+\frac{1}{\theta\gamma\xi}\left(  \chi_{4\lambda}%
\overline{\varphi_{4}}^{^{\prime}}-\chi_{4\lambda}^{^{\prime}}\overline
{\varphi_{4}}\right)  \mid_{h_{3}}^{1}=
\]%
\begin{align*}
&  \chi_{1\lambda}\left(  h_{1},\lambda\right)  \overline{\varphi_{1}%
}^{^{\prime}}\left(  h_{1},\lambda\right)  -\chi_{1\lambda}^{^{\prime}}\left(
h_{1},\lambda\right)  \overline{\varphi_{1}}\left(  h_{1},\lambda\right)
-\chi_{1\lambda}\left(  -1,\lambda\right)  \overline{\varphi_{1}}^{^{\prime}%
}\left(  -1,\lambda\right)  +\\
&  \chi_{1\lambda}^{^{\prime}}\left(  -1,\lambda\right)  \overline{\varphi
_{1}}\left(  -1,\lambda\right)  +\frac{1}{\theta}\left(  \chi_{2\lambda
}\left(  h_{2},\lambda\right)  \overline{\varphi_{2}}^{^{\prime}}\left(
h_{2},\lambda\right)  -\chi_{2\lambda}^{^{\prime}}\left(  h_{2},\lambda
\right)  \overline{\varphi_{2}}\left(  h_{2},\lambda\right)  \right.  \\
&  \left.  -\chi_{2\lambda}\left(  h_{1},\lambda\right)  \overline{\varphi
_{2}}^{^{\prime}}\left(  h_{1},\lambda\right)  +\chi_{2\lambda}^{^{\prime}%
}\left(  h_{1},\lambda\right)  \overline{\varphi_{2}}\left(  h_{1}%
,\lambda\right)  \right)  +\frac{1}{\theta\gamma}\left(  \chi_{3\lambda
}\left(  h_{3},\lambda\right)  \overline{\varphi_{3}}^{^{\prime}}\left(
h_{3},\lambda\right)  \right.  \\
&  \left.  -\chi_{3\lambda}^{^{\prime}}\left(  h_{3},\lambda\right)
\overline{\varphi_{3}}\left(  h_{3},\lambda\right)  -\chi_{3\lambda}\left(
h_{2},\lambda\right)  \overline{\varphi_{3}}^{^{\prime}}\left(  h_{2}%
,\lambda\right)  +\chi_{3\lambda}^{^{\prime}}\left(  h_{2},\lambda\right)
\overline{\varphi_{3}}\left(  h_{2},\lambda\right)  \right)  +\\
&  \frac{1}{\theta\gamma\xi}\left(  \chi_{4\lambda}\left(  1,\lambda\right)
\overline{\varphi_{4}}^{^{\prime}}\left(  1,\lambda\right)  -\chi_{4\lambda
}^{^{\prime}}\left(  1,\lambda\right)  \overline{\varphi_{4}}\left(
1,\lambda\right)  -\chi_{4\lambda}\left(  h_{3},\lambda\right)  \overline
{\varphi_{4}}^{^{\prime}}\left(  h_{3},\lambda\right)  \right.  \\
&  +\chi_{4\lambda}^{^{\prime}}\left(  h_{3},\lambda\right)  \overline
{\varphi_{4}}\left(  h_{3},\lambda\right)  \\
&  =\alpha_{1}\chi_{1\lambda}\left(  -1,\lambda\right)  +\alpha_{2}%
\chi_{1\lambda}^{^{\prime}}\left(  -1,\lambda\right)  +\chi_{1\lambda}\left(
h_{1},\lambda\right)  \overline{\varphi_{1}}^{^{\prime}}\left(  h_{1}%
,\lambda\right)  -\chi_{1\lambda}^{^{\prime}}\left(  h_{1},\lambda\right)
\overline{\varphi_{1}}\left(  h_{1},\lambda\right)  \\
&  +\frac{1}{\theta}\left(  \chi_{2\lambda}\left(  h_{2},\lambda\right)
\overline{\varphi_{2}}^{^{\prime}}\left(  h_{2},\lambda\right)  -\chi
_{2\lambda}^{^{\prime}}\left(  h_{2},\lambda\right)  \overline{\varphi_{2}%
}\left(  h_{2},\lambda\right)  -\chi_{2\lambda}\left(  h_{1},\lambda\right)
\overline{\varphi_{2}}^{^{\prime}}\left(  h_{1},\lambda\right)  +\right.  \\
&  \left.  \chi_{2\lambda}^{^{\prime}}\left(  h_{1},\lambda\right)
\overline{\varphi_{2}}\left(  h_{1},\lambda\right)  \right)  +\frac{1}%
{\theta\gamma}\left(  \chi_{3\lambda}\left(  h_{3},\lambda\right)
\overline{\varphi_{3}}^{^{\prime}}\left(  h_{3},\lambda\right)  -\chi
_{3\lambda}\left(  h_{3},\lambda\right)  \overline{\varphi_{3}}\left(
h_{3},\lambda\right)  \right)  \\
&  -\frac{1}{\theta\gamma}\left(  \chi_{3\lambda}\left(  h_{2},\lambda\right)
\overline{\varphi_{3}}^{^{\prime}}\left(  h_{2},\lambda\right)  -\chi
_{3\lambda}^{^{\prime}}\left(  h_{2},\lambda\right)  \overline{\varphi_{3}%
}\left(  h_{2},\lambda\right)  \right)  +\\
&  \frac{1}{\theta\gamma\xi}\left(  \beta_{2}^{^{\prime}}\overline{\varphi
_{4}}^{^{\prime}}\left(  1,\lambda\right)  -\beta_{1}^{^{\prime}}%
\overline{\varphi_{4}}\left(  1,\lambda\right)  \right)  -\frac{1}%
{\theta\gamma\xi}\left(  \chi_{4\lambda}\left(  h_{3},\lambda\right)
\overline{\varphi_{4}}^{^{\prime}}\left(  h_{3},\lambda\right)  -\chi
_{4\lambda}^{^{\prime}}\left(  h_{3},\lambda\right)  \overline{\varphi_{4}%
}\left(  h_{3},\lambda\right)  \right)  \\
&  =\alpha_{1}\chi_{1\lambda}\left(  -1,\lambda\right)  +\alpha_{2}%
\chi_{1\lambda}^{^{\prime}}\left(  -1,\lambda\right)  +\frac{1}{\theta
\gamma\xi}\left(  \beta_{2}^{^{\prime}}\overline{\varphi_{4}}^{^{\prime}%
}\left(  1,\lambda\right)  -\beta_{1}^{^{\prime}}\overline{\varphi_{4}}\left(
1,\lambda\right)  \right)  .
\end{align*}
Note that%
\[
w^{^{\prime}}\left(  \lambda\right)  =\alpha_{2}\chi_{1\lambda}^{^{\prime}%
}\left(  -1,\lambda\right)  +\alpha_{1}\chi_{1\lambda}\left(  -1,\lambda
\right)
\]
Therefore, $\left(  9\right)  $ becomes%
\begin{equation}
w^{^{\prime}}\left(  \lambda\right)  =\left\langle \chi,\phi\right\rangle
_{1}+2it\left\langle \chi_{\lambda},\phi\right\rangle _{1}-\frac{1}%
{\theta\gamma\xi}\left(  \beta_{2}^{^{\prime}}\overline{\varphi_{4}}%
^{^{\prime}}\left(  1,\lambda\right)  -\beta_{1}^{^{\prime}}\overline
{\varphi_{4}}\left(  1,\lambda\right)  \right)  .\tag{10}%
\end{equation}
Now we consider the simplicity of the eigenvalues of $\left(  1\right)
-\left(  9\right)  $. Let $\mu$ be arbitrary zero of $w\left(  \lambda\right)
$. By Corollary 1.1, $\mu$ is real. S\.{I}nce%
\[
w\left(  \mu\right)  =\left\vert
\begin{array}
[c]{cc}%
\varphi_{1}\left(  x,\mu\right)   & \chi_{1}\left(  x,\mu\right)  \\
\varphi_{1}^{^{\prime}}\left(  x,\mu\right)   & \chi_{1}^{^{\prime}}\left(
x,\mu\right)
\end{array}
\right\vert =0,
\]
we have $\varphi_{1}\left(  x,\mu\right)  =c_{1}\chi_{1}\left(  x,\mu\right)
$ $\left(  c_{1}\neq0\right)  ,$ $\varphi_{2}\left(  x,\mu\right)  =c_{2}%
\chi_{2}\left(  x,\mu\right)  $ $\left(  c_{2}\neq0\right)  ,$ $\varphi
_{3}\left(  x,\mu\right)  =c_{3}\chi_{3}\left(  x,\mu\right)  $ $\left(
c_{3}\neq0\right)  $ and $\varphi_{4}\left(  x,\mu\right)  =c_{4}\chi
_{4}\left(  x,\mu\right)  $ $\left(  c_{4}\neq0\right)  $ where $c_{1}%
,c_{2},c_{3},c_{4}\in%
%TCIMACRO{\U{2102} }%
%BeginExpansion
\mathbb{C}
%EndExpansion
.$ From%
\begin{align*}
\varphi_{2}\left(  h_{1},\mu\right)   &  =c_{1}\left(  \alpha_{3}\chi
_{1}\left(  h_{1},\mu\right)  +\beta_{3}\chi_{1}^{^{\prime}}\left(  h_{1}%
,\mu\right)  \right)  =c_{1}\chi_{2}\left(  h_{1},\mu\right)  ,\\
\varphi_{2}^{^{\prime}}\left(  h_{1},\mu\right)   &  =c_{1}\left(  \alpha
_{4}\chi_{1}\left(  h_{1},\mu\right)  +\beta_{4}\chi_{1}^{^{\prime}}\left(
h_{1},\mu\right)  \right)  =c_{1}\chi_{2}^{^{\prime}}\left(  h_{1},\mu\right)
,\\
\varphi_{3}\left(  h_{2},\mu\right)   &  =c_{2}\left(  \alpha_{5}\chi
_{2}\left(  h_{2},\mu\right)  +\beta_{5}\chi_{2}^{^{\prime}}\left(  h_{2}%
,\mu\right)  \right)  =c_{2}\chi_{3}\left(  h_{2},\mu\right)  ,\\
\varphi_{3}^{^{\prime}}\left(  h_{2},\mu\right)   &  =c_{2}\left(  \alpha
_{6}\chi_{2}\left(  h_{2},\mu\right)  +\beta_{6}\chi_{2}^{^{\prime}}\left(
h_{2},\mu\right)  \right)  =c_{2}\chi_{3}^{^{\prime}}\left(  h_{2},\mu\right)
,\\
\varphi_{4}\left(  h_{3},\mu\right)   &  =c_{3}\left(  \alpha_{7}\chi
_{3}\left(  h_{3},\mu\right)  +\beta_{7}\chi_{3}^{^{\prime}}\left(  h_{3}%
,\mu\right)  \right)  =c_{3}\chi_{4}\left(  h_{3},\mu\right)  ,\\
\varphi_{4}^{^{\prime}}\left(  h_{3},\mu\right)   &  =c_{3}\left(  \alpha
_{8}\chi_{3}\left(  h_{3},\mu\right)  +\beta_{8}\chi_{3}^{^{\prime}}\left(
h_{3},\mu\right)  \right)  =c_{3}\chi_{4}^{^{\prime}}\left(  h_{3},\mu\right)
,
\end{align*}
we get $c_{1}=c_{2}=c_{3}=c_{4}\neq0.$ Thus, simple calculations using
$\left(  10\right)  $ and the initial values of $\chi_{4}$ at $x=1$ give%
\begin{align*}
w^{^{\prime}}\left(  \mu\right)    & =\overline{c_{1}}\left(  \frac{1}%
{p_{1}^{2}}%
%TCIMACRO{\dint \limits_{-1}^{h_{1}}}%
%BeginExpansion
{\displaystyle\int\limits_{-1}^{h_{1}}}
%EndExpansion
\left\vert \chi_{1}\left(  x,\mu\right)  \right\vert ^{2}dx+\frac{1}{p_{2}%
^{2}\theta}%
%TCIMACRO{\dint \limits_{h_{1}}^{h_{2}}}%
%BeginExpansion
{\displaystyle\int\limits_{h_{1}}^{h_{2}}}
%EndExpansion
\left\vert \chi_{2}\left(  x,\mu\right)  \right\vert ^{2}dx+\frac{1}{p_{3}%
^{2}\theta\gamma}%
%TCIMACRO{\dint \limits_{h_{2}}^{h_{3}}}%
%BeginExpansion
{\displaystyle\int\limits_{h_{2}}^{h_{3}}}
%EndExpansion
\left\vert \chi_{3}\left(  x,\mu\right)  \right\vert ^{2}dx\right.  \\
& \left.  +\frac{1}{p_{4}^{2}\theta\gamma\xi}%
%TCIMACRO{\dint \limits_{h_{3}}^{1}}%
%BeginExpansion
{\displaystyle\int\limits_{h_{3}}^{1}}
%EndExpansion
\left\vert \chi_{4}\left(  x,\mu\right)  \right\vert ^{2}dx+\frac{\rho}%
{\theta\gamma\xi}.\right)
\end{align*}
Note that $\rho>0,$ $\theta>0,$ $\gamma>0,$ $\xi>0$ and $c_{1}\neq0,$ so
$w^{^{\prime}}\left(  \mu\right)  \neq0.$ Hence, the analytic multiplicity of
$\mu$ is one. By Lemma 2.1, the proof is completed.
\end{proof}

\textbf{Theorem 2.2 }\textit{All eigenvalues of }$\left(  1\right)  -\left(
9\right)  $\textit{ are geometrically simple.}

\begin{proof}
If $f$ and $g$ are two eigenfunctions for an eigenvalue $\lambda_{\ast}$ of
$\left(  1\right)  -\left(  9\right)  $, then $\left(  2\right)  $ implies
that $f\left(  -1\right)  =Kg\left(  -1\right)  $ and $f^{^{\prime}}\left(
-1\right)  =Kg^{^{\prime}}\left(  -1\right)  $ for some constant $K\in%
%TCIMACRO{\U{211d} }%
%BeginExpansion
\mathbb{R}
%EndExpansion
.$ By the uniqueness \ theorem for solutions of ordinary differential equation
and the transmission conditions $\left(  4\right)  -\left(  9\right)  ,$ we
have that $f=Kg$ on $\left[  -1,1\right]  .$ Thus the geometric multiplicity
of $\lambda_{\ast}$ is one.
\end{proof}

\section{Completeness of eigenfunctions}

\textbf{Theorem 3.1 }\textit{The operator }$A$\textit{ has only point
spectrum, i.e., }$\sigma\left(  A\right)  =$\textit{ }$\sigma_{\rho}\left(
A\right)  .$

\begin{proof}
It suffices to prove that if $\eta$ is not an eigenvalue of $A,$ then $\eta
\in\rho\left(  A\right)  .$ Since $A$ is self-adjoint, we only consider a real
$\eta.$ We investigate the equation $\left(  A-\eta\right)  Y=F\in H,$ where
$F=\left(  f,h\right)  .$

Let us consider the initial-value problem%
\begin{equation}
\left\{
\begin{array}
[c]{c}%
\tau y-\eta y=f,\text{ \ }x\in I,\\
\alpha_{1}y\left(  -1\right)  +\alpha_{2}y^{^{\prime}}\left(  -1\right)  =0,\\
y\left(  h_{1}+0\right)  =\alpha_{3}y\left(  h_{1}-0\right)  +\beta
_{3}y^{^{\prime}}\left(  h_{1}-0\right)  ,\\
y^{^{\prime}}\left(  h_{1}+0\right)  =\alpha_{4}y\left(  h_{1}-0\right)
+\beta_{4}y^{^{\prime}}\left(  h_{1}-0\right)  ,\\
y\left(  h_{2}+0\right)  =\alpha_{5}y\left(  h_{2}-0\right)  +\beta
_{5}y^{^{\prime}}\left(  h_{2}-0\right)  ,\\
y^{^{\prime}}\left(  h_{2}+0\right)  =\alpha_{6}y\left(  h_{2}-0\right)
+\beta_{6}y^{^{\prime}}\left(  h_{2}-0\right)  ,\\
y\left(  h_{3}+0\right)  =\alpha_{7}y\left(  h_{3}-0\right)  +\beta
_{7}y^{^{\prime}}\left(  h_{3}-0\right)  ,\\
y^{^{\prime}}\left(  h_{3}+0\right)  =\alpha_{8}y\left(  h_{3}-0\right)
+\beta_{8}y^{^{\prime}}\left(  h_{3}-0\right)  ,
\end{array}
\right.  \tag{11}%
\end{equation}
Let $u\left(  x\right)  $ be the solution of the equation $\tau u-\eta u=0$
satisfying%
\begin{align*}
u\left(  -1\right)   &  =\alpha_{2},\text{ }u^{^{\prime}}\left(  -1\right)
=-\alpha_{1},\\
u\left(  h_{1}+0\right)   &  =\alpha_{3}u\left(  h_{1}-0\right)  +\beta
_{3}u^{^{\prime}}\left(  h_{1}-0\right)  ,\\
u^{^{\prime}}\left(  h_{1}+0\right)   &  =\alpha_{4}u\left(  h_{1}-0\right)
+\beta_{4}u^{^{\prime}}\left(  h_{1}-0\right)  ,\\
u\left(  h_{2}+0\right)   &  =\alpha_{5}u\left(  h_{2}-0\right)  +\beta
_{5}u^{^{\prime}}\left(  h_{2}-0\right)  ,\\
u^{^{\prime}}\left(  h_{2}+0\right)   &  =\alpha_{6}u\left(  h_{2}-0\right)
+\beta_{6}u^{^{\prime}}\left(  h_{2}-0\right)  ,\\
u\left(  h_{3}+0\right)   &  =\alpha_{7}u\left(  h_{3}-0\right)  +\beta
_{7}u^{^{\prime}}\left(  h_{3}-0\right)  ,\\
u^{^{\prime}}\left(  h_{3}+0\right)   &  =\alpha_{8}u\left(  h_{3}-0\right)
+\beta_{8}u^{^{\prime}}\left(  h_{3}-0\right)  .
\end{align*}
In fact,%
\[
u\left(  x\right)  =\left\{
\begin{array}
[c]{c}%
u_{1}\left(  x\right)  ,\text{ \ }x\in\left[  -1,h_{1}\right)  ,\\
u_{2}\left(  x\right)  ,\text{ \ }x\in\left(  h_{1},h_{2}\right)  ,\\
u_{3}\left(  x\right)  ,\text{ \ }x\in\left(  h_{2},h_{3}\right)  ,\\
u_{4}\left(  x\right)  ,\text{ \ }x\in\left(  h_{3},1\right]  ,
\end{array}
\right.
\]
where $u_{1}\left(  x\right)  $ is the unique solution of the initial-value
problem%
\[
\left\{
\begin{array}
[c]{c}%
-p_{1}^{2}u^{^{\prime\prime}}+q(x)u=\eta u,\text{ \ \ }x\in\left[
-1,h_{1}\right)  ,\\
u\left(  -1\right)  =\alpha_{2},\text{ }u^{^{\prime}}\left(  -1\right)
=-\alpha_{1};
\end{array}
\right.
\]
$u_{2}\left(  x\right)  $ is the unique solution of the problem%
\[
\left\{
\begin{array}
[c]{c}%
-p_{2}^{2}u^{^{\prime\prime}}+q(x)u=\eta u,\text{ \ \ }x\in\left(  h_{1}%
,h_{2}\right)  ,\\
u_{2}\left(  h_{1}\right)  =\alpha_{3}u_{1}\left(  h_{1}\right)  +\beta
_{3}u_{1}^{^{\prime}}\left(  h_{1}\right)  ,\\
u_{2}^{^{\prime}}\left(  h_{1}\right)  =\alpha_{4}u_{1}\left(  h_{1}\right)
+\beta_{4}u_{1}^{^{\prime}}\left(  h_{1}\right)  ;
\end{array}
\right.
\]
$u_{3}\left(  x\right)  $ is the unique solution of the problem%
\[
\left\{
\begin{array}
[c]{c}%
-p_{3}^{2}u^{^{\prime\prime}}+q(x)u=\eta u,\text{ \ \ }x\in\left(  h_{2}%
,h_{3}\right)  ,\\
u_{3}\left(  h_{2}\right)  =\alpha_{5}u_{2}\left(  h_{2}\right)  +\beta
_{5}u_{2}^{^{\prime}}\left(  h_{2}\right)  ,\\
u_{3}^{^{\prime}}\left(  h_{2}\right)  =\alpha_{6}u_{2}\left(  h_{2}\right)
+\beta_{6}u_{2}^{^{\prime}}\left(  h_{2}\right)  .
\end{array}
\right.
\]
and $u_{4}\left(  x\right)  $ is the unique solution of the problem%
\[
\left\{
\begin{array}
[c]{c}%
-p_{4}^{2}u^{^{\prime\prime}}+q(x)u=\eta u,\text{ \ \ }x\in\left(
h_{3},1\right]  ,\\
u_{3}\left(  h_{3}\right)  =\alpha_{7}u_{3}\left(  h_{3}\right)  +\beta
_{7}u_{3}^{^{\prime}}\left(  h_{3}\right)  ,\\
u_{3}^{^{\prime}}\left(  h_{3}\right)  =\alpha_{8}u_{3}\left(  h_{3}\right)
+\beta_{8}u_{3}^{^{\prime}}\left(  h_{3}\right)  .
\end{array}
\right.
\]
Let%
\[
w\left(  x\right)  =\left\{
\begin{array}
[c]{c}%
w_{1}\left(  x\right)  ,\text{ \ }x\in\left[  -1,h_{1}\right)  ,\\
w_{2}\left(  x\right)  ,\text{ \ }x\in\left(  h_{1},h_{2}\right)  ,\\
w_{3}\left(  x\right)  ,\text{ \ }x\in\left(  h_{2},h_{3}\right)  ,\\
w_{4}\left(  x\right)  ,\text{ \ }x\in\left(  h_{3},1\right]  ,
\end{array}
\right.
\]
be a solution of $\tau w-\eta w=f$ satisfying%
\[%
\begin{array}
[c]{c}%
\alpha_{1}w\left(  -1\right)  +\alpha_{2}w^{^{\prime}}\left(  -1\right)  =0,\\
w\left(  h_{1}+0\right)  =\alpha_{3}w\left(  h_{1}-0\right)  +\beta
_{3}w^{^{\prime}}\left(  h_{1}-0\right)  ,\\
w^{^{\prime}}\left(  h_{1}+0\right)  =\alpha_{4}w\left(  h_{1}-0\right)
+\beta_{4}w^{^{\prime}}\left(  h_{1}-0\right)  ,\\
w\left(  h_{2}+0\right)  =\alpha_{5}w\left(  h_{2}-0\right)  +\beta
_{5}w^{^{\prime}}\left(  h_{2}-0\right)  ,\\
w^{^{\prime}}\left(  h_{2}+0\right)  =\alpha_{6}w\left(  h_{2}-0\right)
+\beta_{6}w^{^{\prime}}\left(  h_{2}-0\right)  ,\\
w\left(  h_{3}+0\right)  =\alpha_{7}w\left(  h_{3}-0\right)  +\beta
_{7}w^{^{\prime}}\left(  h_{3}-0\right)  ,\\
w^{^{\prime}}\left(  h_{3}+0\right)  =\alpha_{8}w\left(  h_{3}-0\right)
+\beta_{8}w^{^{\prime}}\left(  h_{3}-0\right)  .
\end{array}
\]
Then, $\left(  11\right)  $ has the general solution%
\begin{equation}
y\left(  x\right)  =\left\{
\begin{array}
[c]{c}%
du_{1}+w_{1},\text{ \ }x\in\left[  -1,h_{1}\right)  ,\\
du_{2}+w_{2},\text{ \ }x\in\left(  h_{1},h_{2}\right)  ,\\
du_{3}+w_{3},\text{ \ }x\in\left(  h_{2},h_{3}\right)  ,\\
du_{4}+w_{4},\text{ \ }x\in\left(  h_{3},1\right]  ,
\end{array}
\right.  \tag{12}%
\end{equation}
where $d\in%
%TCIMACRO{\U{2102} }%
%BeginExpansion
\mathbb{C}
%EndExpansion
.$

Since $\eta$ is not an eigenvalue of the problem $\left(  1\right)  -\left(
7\right)  ,$ we have%
\begin{equation}
\eta\left[  \beta_{1}^{^{\prime}}u_{2}\left(  1\right)  -\beta_{2}^{^{\prime}%
}u_{2}^{^{\prime}}(1)\right]  +\left[  \beta_{1}u_{2}\left(  1\right)
-\beta_{2}u_{2}^{^{\prime}}(1)\right]  \neq0. \tag{13}%
\end{equation}
The second component of $\left(  A-\eta\right)  Y=F$ involves the equation%
\[
-N\left(  y\right)  -\eta N^{^{\prime}}\left(  y\right)  =h,
\]
namely,%
\begin{equation}
\left[  -\beta_{1}y\left(  1\right)  +\beta_{2}y^{^{\prime}}(1)\right]
-\eta\left[  \beta_{1}^{^{\prime}}y\left(  1\right)  -\beta_{2}^{^{\prime}%
}y^{^{\prime}}(1)\right]  =h. \tag{14}%
\end{equation}
Substituting $\left(  12\right)  $ into $\left(  14\right)  $, we get%
\begin{align*}
&  \left(  \beta_{2}u_{2}^{^{\prime}}\left(  1\right)  -\beta_{1}u_{2}\left(
1\right)  +\eta\beta_{2}^{^{\prime}}u_{2}^{^{\prime}}\left(  1\right)
-\eta\beta_{1}^{^{\prime}}u_{2}\left(  1\right)  \right)  d\\
&  =h+\beta_{1}w_{2}\left(  1\right)  -\beta_{2}w_{2}^{^{\prime}}\left(
1\right)  +\eta\beta_{1}^{^{\prime}}w_{2}\left(  1\right)  -\eta\beta
_{2}^{^{\prime}}w_{2}^{^{\prime}}\left(  1\right)
\end{align*}
In view of $\left(  13\right)  $, we know that $d$ is uniquely solvable.
Therefore, $y$ is uniquely determined.

The above arguments show that $\left(  A-\eta I\right)  ^{-1}$ is defined on
all of $H,$ where $I$ is identity matrix. We obtain that $\left(  A-\eta
I\right)  ^{-1}$ is bounded by Theorem 1.1 and the Closed Graph Theorem. Thus,
$\eta\in\rho\left(  A\right)  .$ Therefore, $\sigma\left(  A\right)
=\sigma_{\rho}\left(  A\right)  .$
\end{proof}

\textbf{Lemma 3.1 }$\left[  10\right]  $\textit{ The eigenvalues of the
boundary value problem }$\left(  1\right)  -\left(  9\right)  $\textit{ are
bounded below, and they are countably infinite and can cluster only at
}$\infty.$

For every $\delta\in%
%TCIMACRO{\U{211d} }%
%BeginExpansion
\mathbb{R}
%EndExpansion
\setminus\sigma_{\rho}\left(  A\right)  ,$ we have the following immediate conclusion.

\textbf{Lemma 3.2 }\textit{Let }$\lambda$\textit{ be an eigenvalue of
}$A-\delta I,$\textit{ and }$V$\textit{ a corresponding eigenfunction. Then,
}$\lambda^{-1}$\textit{ is an eigenvalue of }$\left(  A-\delta I\right)
^{-1},$\textit{ and }$V$\textit{ is a corresponding eigenfunction. The
converse is also true.}

On the other hand, if $\mu$ is an eigenvalue of $A$ and $U$ is a corresponding
eigenfunction, then $\mu-\delta$ is an eigenvalue of $A-\delta I,$ and $U$ is
a corresponding eigenfunction. The converse is also true. Accordingly, the
discussion about the completeness of the eigenfunctions of $A$ is equivalent
to considering the corresponding property of $\left(  A-\delta I\right)
^{-1}$.

By Lemma 1.1, Lemma 3.1 and Corollary 1.1, we suppose that $\left\{
\lambda_{n};\text{ }n\in%
%TCIMACRO{\U{2115} }%
%BeginExpansion
\mathbb{N}
%EndExpansion
\right\}  $ is the real sequence of eigenvalues of $A,$ then $\left\{
\lambda_{n}-\delta;\text{ }n\in%
%TCIMACRO{\U{2115} }%
%BeginExpansion
\mathbb{N}
%EndExpansion
\right\}  $ is the sequence of eigenvalues of $A-\delta I$. We may assume that%
\[
\left\vert \lambda_{1}-\delta\right\vert \leq\left\vert \lambda_{2}%
-\delta\right\vert \leq...\leq\left\vert \lambda_{n}-\delta\right\vert
\leq...\rightarrow\infty.
\]
Let $\left\{  \mu_{n};\text{ }n\in%
%TCIMACRO{\U{2115} }%
%BeginExpansion
\mathbb{N}
%EndExpansion
\right\}  $ be the sequence of eigenvalues of $\left(  A-\delta I\right)
^{-1}.$ Then $\mu_{n}=\left(  \lambda_{n}-\delta\right)  ^{-1}$ and%
\[
\left\vert \mu_{1}\right\vert \geq\left\vert \mu_{2}\right\vert \geq
...\geq\left\vert \mu_{n}\right\vert \geq...\rightarrow0.
\]
Note that $0$ is not an eigenvalue of $\left(  A-\delta I\right)  ^{-1}.$

\textbf{Theorem 3.2 }\textit{The operator }$A$\textit{ has compact resolvents,
i.e, for each }$\delta\in%
%TCIMACRO{\U{211d} }%
%BeginExpansion
\mathbb{R}
%EndExpansion
\setminus\sigma_{\rho}\left(  A\right)  ,$\textit{ }$\left(  A-\delta
I\right)  ^{-1}$\textit{ is compact on }$H.$

\begin{proof}
Let $\left\{  \mu_{1},\mu_{2},...\right\}  $ be the eigenvalues of $\left(
A-\delta I\right)  ^{-1},$ and $\left\{  P_{1},P_{2},...\right\}  $ the
orthogonal projections of finite rank onto the corresponding eigenspaces.
Since $\left\{  \mu_{1},\mu_{2},...\right\}  $ is a bounded sequence and all
$P_{n}\prime s$ are mutually orthogonal, we have $\sum_{n=1}^{\infty}\mu
_{n}P_{n}$ is strongly convergent to the bounded operator $\left(  A-\delta
I\right)  ^{-1},$ i.e., $\left(  A-\delta I\right)  ^{-1}=\sum_{n=1}^{\infty
}\mu_{n}P_{n}.$ Because for every $\alpha>0,$ the number of $\mu_{n}\prime s$
satisfying $\left\vert \mu_{n}\right\vert >\alpha$ is finite, and all
$P_{n}\prime s$ are of finite rank, we obtain that $\left(  A-\delta I\right)
^{-1}$ is compact.
\end{proof}

In terms of the above statements and the spectral theorem for compact
operators, we obtain the following theorem.

\textbf{Theorem 3.3 }\textit{The eigenfunctions of the problem }$\left(
1\right)  -\left(  9\right)  ,$\textit{ augmented to become eigenfunctions of
}$A,$\textit{ are complete in }$H$\textit{, i.e., if we let}%
\[
\left\{  \Phi_{n}=\left(  \phi_{n}\left(  x\right)  ,N^{^{\prime}}\left(
\phi_{n}\right)  \right)  ;\text{ }n\in%
%TCIMACRO{\U{2115} }%
%BeginExpansion
\mathbb{N}
%EndExpansion
\right\}
\]
\textit{ be a maximum set of orthonormal eigenfunctions of }$A,$\textit{ where
}$\left\{  \phi_{n}\left(  x\right)  ;\text{ }n\in%
%TCIMACRO{\U{2115} }%
%BeginExpansion
\mathbb{N}
%EndExpansion
\right\}  $\textit{ are eigenfunctions of }$\left(  1\right)  -\left(
9\right)  ,$\textit{ then for all }$F\in H,$\textit{ }$F=\sum_{n=1}^{\infty
}\left\langle F,\Phi_{n}\right\rangle \Phi_{n}.$


\begin{thebibliography}{99}                                                                                               %


\bibitem {}{\small C. T. Fulton, Two-point boundary value problems with
eigenvalue parameter contained in the boundary conditions. Proc. Royal Soc.
Edinburgh 77A (1977), 293-308.}

\bibitem {}{\small A. V. Likov, Y. A. Mikhalilov, The Theory of Heat and Mass
Transfer, Qosenergaizdat, 1963. (In Russian).}

\bibitem {}{\small E. \c{S}en, A. Bayramov, Calculation of eigenvalues and
eigenfunctions of a discontinuous boundary value problem with retarded
argument which contains a spectral parameter in the boundary condition,
Mathematical and Computer Modelling, 54 (2011) 3090-3097.}

\bibitem {}{\small E. \c{S}en, A. Bayramov, On calculation of eigenvalues and
eigenfunctions of a Sturm-Liouville type problem with retarded argument which
contains a spectral parameter in the boundary condition, Journal of
Inequalities and Applications, Vol. 2011 (1) (2011) 1-9.}

\bibitem {}{\small P. A. Binding, P. J. Browne, B. A. Watson, Sturm-Liouville
problems with boundary conditions rationally dependent on the eigenparameter,
II, J. Comput. Appl. Math. 148 (2002), 147-168.}

\bibitem {}{\small P. A. Binding, Patrick. J. Browne, Oscillation theory for
indefinite Sturm-Liouville problems with eigenparameter dependent boundary
conditions. Proc. Royal Soc. Edinburg 127A (1997), 1123-1136.}

\bibitem {}{\small P. A. Binding, R. Hryniv, H. Langer, Elliptic eigenvalue
problems with eigenparameter dependent boundary conditions. J. Differential
Equations 174 (2001), 30-54.}

\bibitem {}{\small M. Demirci, Z. Akdo\u{g}an, O. Sh. Mukhtarov, Asymptotic
behavior of eigenvalues and eigenfunctions of one discontinuous boundary-value
problem. International J. Computational Cognition 2(3) (2004), 101-113.}

\bibitem {}{\small D. B. Hinton, An expansion theorem for an eigenvalue
problem with eigenvalue parameter in the boundary condition. Quart. J. Math.
(Oxford) 30 (1979), 33-42.}

\bibitem {}{\small M. Kadakal, O.Sh. Mukhtarov, Sturm-Liouville problems with
discontinuities at two points, Comput. Math. Appl.,54 (2007) 1367-1379.}

\bibitem {}{\small O. Sh. Mukhtarov, S. Yakubov, Problems for differential
equations with transmission conditions, Applicable Anal. 81 (2002),
1033-1064.}

\bibitem {}{\small I. Titeux, Y. Yakubov, Completeness of root functions for
thermal conduction in a strip with piecewise continuous coefficients, Math.
Models Methods Appl. Sc. 7 (1997), 1035-1050.}

\bibitem {}{\small O. Sh. Mukhtarov, M. Kadakal, Spectral properties of one
Sturm-Liouville type problem with discontinuous weight, (Russian) Sibirsk.
Math. Zh., 46 (4) (2005), 681-694.}

\bibitem {}{\small Y. S. Li, J. Sun, Z. Wang, On the complete description of
self-adjoint boundary conditions of the Schrodinger operator with a }%
$\delta\left(  x\right)  ${\small or }$\delta^{^{\prime}}\left(  x\right)
${\small interaction, Symposium of the Fifth Conference of Mathematics Society
of Inner Mongolia, Inner Mongolia Univ. Press, 1995, 27-30.}

\bibitem {}{\small A. Zettl, Adjoint and self-adjoint boundary value problems
with interface conditions, SIAM J. Appl. Math. 16 (1968), 851-859.}

\bibitem {}{\small F. Gesztesy, W. Kirsch, One-dimensional Schrodinger
operators with interactions on a discerete set, J. Reine Angew. Math. 362
(1985), 28-50.}

\bibitem {}{\small A. Wang, J. Sun, P. Gao, Completeness of Eigenfunctions of
Sturm-Liouville Problems with Transmission Conditions, J. Spectral Math. Appl.
(2006)}

\bibitem {}{\small A. Wang, J. Sun, X. Hao, S. Yao, Completeness of
Eigenfunctions of Sturm-Liouville Problems with Transmission Conditions,
Methods and Application of Analysis 16 (3) (2009), 299-312.}

\bibitem {}{\small J. Weidmann, Spectral Theory of Ordinary Differential
Operators, Lecture Notes in Math. 1258, Springer-Verlag, Berline, 1987.}

\bibitem {}{\small M. A. Naimark, Linear Differential Operators, part II.
Harrap, London, 1968.}
\end{thebibliography}
\end{document}